# MULTISCALE LOCAL CHANGE POINT DETECTION WITH APPLICATIONS TO VALUE-AT-RISK

By Vladimir Spokoiny

*Weierstrass-Institute and Humboldt University Berlin*

This paper offers a new approach to modeling and forecasting of nonstationary time series with applications to volatility modeling for financial data. The approach is based on the assumption of local homogeneity: for every time point, there exists a historical *interval of homogeneity*, in which the volatility parameter can be well approximated by a constant. The proposed procedure recovers this interval from the data using the local change point (LCP) analysis. Afterward, the estimate of the volatility can be simply obtained by local averaging. The approach carefully addresses the question of choosing the tuning parameters of the procedure using the so-called "propagation" condition. The main result claims a new "oracle" inequality in terms of the modeling bias which measures the quality of the local constant approximation. This result yields the optimal rate of estimation for smooth and piecewise constant volatility functions. Then, the new procedure is applied to some data sets and a comparison with a standard GARCH model is also provided. Finally, we discuss applications of the new method to the Value at Risk problem. The numerical results demonstrate a very reasonable performance of the new method.

**1. Introduction.** This paper presents a novel approach to modeling of nonstationary time series based on the local parametric assumption, which means that the underlying process having an arbitrary nonstationary structure can, however, be well approximated by a simple time-homogeneous parametric time series within some time interval.

Since the seminal papers of Engle (1982) and Bollerslev (1986), modeling the dynamic features of the variance of financial time series has become one of the most active fields of research in econometrics. New models, different applications and extensions have been proposed, as can be seen by consulting, for example, the monographs of Engle (1995) and of Gouriéroux (1997).

The main idea behind this strain of research is that the volatility clustering effect that is displayed by stock or exchange rate returns can be modeled globally by a stationary process. This approach is somehow restrictive and does not fit some characteristics of the data, in particular the fact that the volatility process appears to be "almost integrated," as it can be seen by usual estimation results and by the very slow decay of the autocorrelations of squared returns. Other global parametric approaches have been proposed by Engle and Bollerslev (1986) and Baillie, Bollerslev and Mikkelsen (1996) in order to include these features in the model. Furthermore, continuous time models, and in particular diffusions and jump diffusions, have also been considered [see, e.g., Andersen, Benzoni and Lund (2002) and Duffie, Pan and Singleton (2000)].

However, Mikosch and Starica (2000b) showed that long memory effects of financial time series can be artificially generated by structural breaks in the parameters. This motivates another modeling approach, which borrows its philosophy mainly from the nonparametric statistics. The main idea consists of using a simple parametric model for describing the conditional distribution of the returns but allowing the parameters of this distribution to be time dependent. The basic assumption of local time homogeneity is that the variability in returns is much higher than the variability in the underlying parameter which allows for estimating this parameter from the most recent historical data. Some examples of this approach can be found in Fan and Gu (2003), Dahlhaus and Rao (2006) and Cheng, Fan and Spokoiny (2003). Furthermore, Mercurio and Spokoiny (2004) proposed a new local adaptive volatility estimation (LAVE) of the unknown volatility from the conditionally heteroskedastic returns. The method is based on pointwise data-driven selection of the interval of homogeneity for every time point. The numerical results demonstrate a reasonable performance of the new method. In particular, it usually outperforms the standard GARCH$(1,1)$ approach. Härdle, Herwartz and Spokoiny (2003) extend this method to estimating the volatility matrix of the multiple returns, and [Mercurio and Torricelli (2003)] apply the same idea in the context of a regression problem.

The aim of the present paper is to develop another approach which, however, applies a similar idea of pointwise adaptive choice of the interval of homogeneity. One essential difference between the LAVE approach from Mercurio and Spokoiny (2004) and the new procedure is in the way of testing the homogeneity of the interval candidate. In this paper, we follow [Grama and Spokoiny (2008)] and systematically apply the approach based on the local multiscale change point analysis. This means that for every historical time point, we test on a structural change at this point for the corresponding scale. The largest interval not containing any change is used for estimation of the parameters of the return distribution. This approach has a number of important advantages of being easy to implement and very sensitive to the structural changes in the return process. We carefully address the question of selecting the tuning parameters of the procedure, which is extremely important for practical applications. The proposed "propagation" approach suggests to tune the parameters under the simple time-homogeneous situation to provide the prescribed performance of the procedure. This way is justified by the theoretical results from Section 4, which claim the "oracle" properties of the resulting estimate in the general situation. Another important feature of the proposed procedure is that it can be easily extended to multiple volatility modeling [cf. Härdle, Herwartz and Spokoiny (2003)].

The change point detection problem for financial time series was considered in Mikosch and Starica (2000a), but they focused on asymptotical properties of the test if only one change point is present. Kitagawa (1987) applied non-Gaussian random walk modeling with heavy tails as the prior for the piecewise constant mean for one-step-ahead prediction of nonstationary time series. However, the mentioned modeling approaches require some essential amount of prior information about the frequency of change points and their size. The new approach proposed in this article does not assume smooth or piecewise constant structure of the underlying process and does not require any prior information. The procedure proposed below in Section 3 focuses on adaptive choice of the interval of homogeneity that allows to proceed in a unified way with smoothly varying coefficient models and change point models.

The proposed LCP approach is quite general and can be applied to many different problems. Grama and Spokoiny (2008) studied the problem of Pareto tail estimation, Giacomini, Härdle and Spokoiny (2007) considered time varying copulae estimation, Čížek, Härdle and Spokoiny (2007) applied it to compare the performance of global and time varying ARCH and GARCH specifications. A comprehensive study of the general LCP procedure is to be given in the forthcoming monograph [Spokoiny (2008)].

The theoretical study given in Sections 2 and 4 focuses on two important features of the proposed procedure: stability in the homogeneous situation and sensitivity to spontaneous changes of the model parameter(s). We particularly show that the procedure provides the optimal sensitivity to changes for the prescribed "false alarm" probability. Note that the classical asymptotic methods for stationary time series do not apply in the considered nonstationary situation with possibly small samples required to develop new approaches and tools. Our way of analysis is based on the so-called "small modeling bias" condition, which generalizes the famous bias-variance trade-off. The main result in Theorem 4.7 claims that the procedure delivers the estimation accuracy corresponding to the optimal choice of the historical interval. It is worth mentioning that the result applies to every volatility process, including piecewise constant, smooth varying or mixed structures.

The paper is organized as follows. Section 2 describes the local parametric approach for the volatility modeling and presents some results about the accuracy of the local constant volatility estimation. Section 3 introduces the adaptive modeling procedure. Theoretical properties of the procedure are discussed in the general situation and for two particular cases: a change point model with piecewise constant volatility and a volatility function smoothly varying in time in Section 4. Section 5 illustrates the performances of the new methodology by means of some simulated examples and real data sets. Note that the same procedure with the default setting is applied for all the examples and applications, and it precisely follows the theoretical description. Section 5.4 discusses applications of the new method to the Value at Risk problem.

**2. Volatility modeling. Local parametric approach.** Let $S_t$ be an observed asset process in discrete time, $t = 1, 2, \ldots$, while $R_t$ defines the corresponding return process: $R_t = \log(S_t/S_{t-1})$. We model this process via the *conditional heteroskedasticity* assumption:

$$R_t = \sigma_t \varepsilon_t, \tag{2.1}$$

where $\varepsilon_t$, $t \geq 1$ is a sequence of independent standard Gaussian random variables and $\sigma_t$ is the *volatility* process which is in general a predictable random process, that is, $\sigma_t$ is measurable with respect to $\mathcal{F}_{t-1}$ with $\mathcal{F}_{t-1} = \sigma(R_1, \ldots, R_{t-1})$ ($\sigma$-field generated by the first $t-1$ observations).

In this paper [similar to Mercurio and Spokoiny (2004)] we focus on the problem of filtering the parameter $f(t) = \sigma_t^2$ from the past observations $R_1, \ldots, R_{t-1}$. This problem naturally arises as an important building block for many tasks of financial engineering like Value at Risk or Portfolio Optimization.

We start the theoretical analysis from the simplest homogeneous case, applying the classical maximum likelihood approach. In particular, we show that the corresponding MLE has nice nonasymptotic properties. Later, we indicate how one can extend these nice results to the general nonhomogeneous situation.

2.1. *Parametric modeling.* A *time-homogeneous* (*time-homoskedastic*) model means that $\sigma_t$ is a constant. The process $S_t$ is then a Geometric Brownian motion observed at discrete time moments. For the homogeneous model $\sigma_t^2 \equiv \theta$ with $t \in I$, the squared returns $Y_t = R_t^2$ follow the equation $Y_t = \theta \varepsilon_t^2$, and the parameter $\theta$ can be estimated using the maximum likelihood method

$$\tilde{\theta}_I = \arg\max_{\theta \geq 0} L_I(\theta) = \arg\max_{\theta \geq 0} \sum_{t \in I} \ell(Y_t, \theta),$$

where $\ell(y,\theta) = -(1/2)\log(2\pi\theta) - y/(2\theta)$ is the log-density of the normal distribution with the parameters $(0,\theta)$. This yields

$$(2.2) \qquad L_I(\theta) = -(N_I/2)\log(2\pi\theta) - S_I/(2\theta),$$

where $N_I$ denotes the number of time points in $I$ and $S_I = \sum_{t\in I} Y_t$.

The volatility model is a particular case of an exponential family, so that a closed form representation for the maximum likelihood estimate $\tilde{\theta}_I$ and for the corresponding fitted log-likelihood $L_I(\tilde{\theta}_I)$ are available [see Polzehl and Spokoiny (2006) for more details].

THEOREM 2.1. *For every interval $I$,*
$$\tilde{\theta}_I = S_I/N_I = N_I^{-1}\sum_{t\in I} Y_t.$$

*Moreover, for every $\theta > 0$, the fitted likelihood ratio $L_I(\tilde{\theta},\theta) = \max_{\theta'} L_I(\theta',\theta)$, with $L_I(\theta',\theta) = L_I(\theta') - L_I(\theta)$, satisfies*

$$(2.3) \qquad L_I(\tilde{\theta}_I,\theta) = N_I \mathcal{K}(\tilde{\theta}_I,\theta),$$

*where*

$$\mathcal{K}(\theta',\theta) = -\{\log(\theta'/\theta) + 1 - \theta'/\theta\}/2$$

*is the Kullback–Leibler information for the two normal distributions with variances $\theta'$ and $\theta$.*

PROOF. Both results follow by simple algebra from (2.2). □

REMARK 2.1. The assumption of normality for the innovations $\varepsilon_t$ is often criticized in the financial literature. Our empirical examples in Section 5.2 below also indicate that the tails of estimated innovations are heavier than the normality would imply. However, the estimate $\tilde{\theta}_I$ remains meaningful even for the nonnormal innovations, it is just a quasi-likelihood approach.

THEOREM 2.2 [Polzehl and Spokoiny (2006)]. *Let $f(t) = \theta^*$ for $t \in I$. If the innovations $\varepsilon_t$ are i.i.d. standard normal, then, for any $\mathfrak{z} > 0$,*
$$\mathbb{P}_{\theta^*}(L_I(\tilde{\theta}_I,\theta^*) > \mathfrak{z}) = \mathbb{P}_{\theta^*}(N_I \mathcal{K}(\tilde{\theta}_I,\theta^*) > \mathfrak{z}) \le 2e^{-\mathfrak{z}}.$$

The result can be extended to the case of nonnormal innovations $\varepsilon_t$ under the condition of bounded exponential moments for $\varepsilon_t^2$. The general case can be reduced to this one by some data transformation [see Chen and Spokoiny (2007) for details].

The Kullback–Leibler divergence $\mathcal{K}$ fulfills $\mathcal{K}(\theta',\theta^*) \le I^*|\theta' - \theta^*|^2$ for any point $\theta'$ in a neighborhood of $\theta^*$, where $I^*$ is the maximum of the Fisher information over this neighborhood. Therefore, the result of Theorem 2.2 guarantees that $|\tilde{\theta}_I - \theta^*| \le CN_I^{-1/2}$ with a high probability. Theorem 2.2 can be used for constructing the confidence intervals for the parameter $\theta^*$.

THEOREM 2.3. *If $\mathfrak{z}_\alpha$ satisfies $2e^{-\mathfrak{z}_\alpha} \leq \alpha$, then*
$$\mathcal{E}_\alpha = \{\theta : N_I \mathcal{K}(\tilde{\theta}_I, \theta) \leq \mathfrak{z}_\alpha\}$$
*is an $\alpha$-confidence set for the parameter $\theta^*$.*

Theorem 2.2 claims that the estimation loss measured by $\mathcal{K}(\tilde{\theta}_I, \theta^*)$ is with high probability bounded by $\mathfrak{z}/N_I$, provided that $\mathfrak{z}$ is sufficiently large. Similarly, one can establish a risk bound for a power loss function.

THEOREM 2.4. *Let $R_t$ be i.i.d. from $\mathcal{N}(0, \theta^*)$. Then, for any $r > 0$ and any interval $I$,*
$$\mathbb{E}_{\theta^*} |L_I(\tilde{\theta}_I, \theta^*)|^r \equiv \mathbb{E}_{\theta^*} |N_I \mathcal{K}(\tilde{\theta}_I, \theta^*)|^r \leq \mathfrak{r}_r,$$
*where $\mathfrak{r}_r = 2r \int_{\mathfrak{z} \geq 0} \mathfrak{z}^{r-1} e^{-\mathfrak{z}} d\mathfrak{z} = 2r\Gamma(r)$. Moreover, there exists a constant $C_r$ depending on $r$ only such that for any $\mathfrak{z} \geq 1$ and any other interval $\mathcal{I}$*
$$\mathbb{E}_{\theta^*} |L_\mathcal{I}(\tilde{\theta}_\mathcal{I}, \theta^*)|^r \mathbf{1}(L_I(\tilde{\theta}_I, \theta^*) \geq \mathfrak{z}) \leq C_r \mathfrak{z}^r e^{-\mathfrak{z}}.$$

PROOF. Proof by Theorem 2.2
$$\mathbb{E}_{\theta^*} |L_I(\tilde{\theta}_I, \theta^*)|^r \leq -\int_{\mathfrak{z} \geq 0} \mathfrak{z}^r d\mathbb{P}_{\theta^*}(L(\tilde{\theta}_I, \theta^*) > \mathfrak{z})$$
$$\leq r \int_{\mathfrak{z} \geq 0} \mathfrak{z}^{r-1} \mathbb{P}_{\theta^*}(L_I(\tilde{\theta}_I, \theta^*) > \mathfrak{z}) d\mathfrak{z} \leq 2r \int_{\mathfrak{z} \geq 0} \mathfrak{z}^{r-1} e^{-\mathfrak{z}} d\mathfrak{z}$$
and the first assertion is fulfilled. Similarly, one can show that
$$\mathbb{E}_{\theta^*} |L_I(\tilde{\theta}_I, \theta^*)|^r \mathbf{1}(L_I(\tilde{\theta}_I, \theta^*) \geq \mathfrak{z}) \leq C_r e^{-\mathfrak{z}},$$
where $C_r$ depends on $r$ only. It remains to note that
$$|L_\mathcal{I}(\tilde{\theta}_\mathcal{I}, \theta^*)|^r \mathbf{1}(L_I(\tilde{\theta}_I, \theta^*) \geq \mathfrak{z})$$
$$\leq |L_I(\tilde{\theta}_I, \theta^*)|^r \mathbf{1}(L_I(\tilde{\theta}_I, \theta^*) \geq \mathfrak{z}) + |L_\mathcal{I}(\tilde{\theta}_\mathcal{I}, \theta^*)|^r \mathbf{1}(L_\mathcal{I}(\tilde{\theta}_\mathcal{I}, \theta^*) \geq \mathfrak{z}). \quad \square$$

2.2. *Risk of estimation in nonparametric situation. "Small modeling bias" condition.* This section extends the bound of Theorem 2.4 to the nonparametric model $R_t^2 = f(t)\varepsilon_t^2$ when the function $f(\cdot)$ is not any longer constant even in a vicinity of the reference point $t^\diamond$. We, however, suppose that the function $f(\cdot)$ can be well approximated by a constant $\theta$ at all points $t \in I$.

Let $Z_\theta = d\mathbb{P}/d\mathbb{P}_\theta$ be the likelihood ratio of the underlying measure $\mathbb{P}$ with regard to the parametric measure $\mathbb{P}_\theta$ corresponding to the constant parameter $f(\cdot) \equiv \theta$. Then,
$$\log Z_\theta = \sum_t \log \frac{p(Y_t, f(t))}{p(Y_t, \theta)}.$$

If we restrict our analysis to an interval $I$ and denote by $\mathbb{P}_I$, respectively, $\mathbb{P}_{I,\theta}$, the measure corresponding to the observations $Y_t$ for $t \in I$, then in a similar way

$$\log Z_{I,\theta} := \log \frac{d\mathbb{P}_I}{d\mathbb{P}_{I,\theta}} = \sum_{t \in I} \log \frac{p(Y_t, f(t))}{p(Y_t, \theta)}.$$

To measure the quality of the approximation of the underlying measure $\mathbb{P}_I$ by the parametric measure $\mathbb{P}_{I,\theta}$, define

(2.4) $$\Delta_I(\theta) = \sum_{t \in I} \mathcal{K}(f(t), \theta),$$

where $\mathcal{K}(f(t), \theta)$ means the Kullback–Leibler distance between two parameter values $f(t)$ and $\theta$.

Let $\varrho(\tilde{\theta}_I, \theta)$ be a loss function for an estimate $\tilde{\theta}_I$ constructed from the observations $Y_t$, for $t \in I$. Define also the corresponding risk under the parametric measure $\mathbb{P}_\theta$:

$$\mathcal{R}(\tilde{\theta}_I, \theta) = \mathbb{E}_\theta \varrho(\tilde{\theta}_I, \theta).$$

The next result explains how the risk bounds can be translated from the parametric to the nonparametric situations.

THEOREM 2.5. *Let, for some $\theta \in \Theta$ and some $\Delta \geq 0$,*

(2.5) $$\mathbb{E}\Delta_I(\theta) \leq \Delta.$$

*Then, it holds for any estimate $\tilde{\theta}$ measurable with regard to $\mathcal{F}_I$*

$$\mathbb{E} \log(1 + \varrho(\tilde{\theta}, \theta)/\mathcal{R}(\tilde{\theta}, \theta)) \leq \Delta + 1.$$

PROOF. The proof is based on the following general result.

LEMMA 2.6. *Let $\mathbb{P}$ and $\mathbb{P}_0$ be two measures such that the Kullback–Leibler divergence $\mathbb{E}\log(d\mathbb{P}/d\mathbb{P}_0)$ satisfies*

$$\mathbb{E}\log(d\mathbb{P}/d\mathbb{P}_0) \leq \Delta < \infty.$$

*Then, for any random variable $\zeta$ with $\mathbb{E}_0\zeta < \infty$,*

$$\mathbb{E}\log(1 + \zeta) \leq \Delta + \mathbb{E}_0\zeta.$$

PROOF. By simple algebra one can check that, for any fixed $y$, the maximum of the function $f(x) = xy - x \log x + x$ is attained at $x = e^y$, leading to the inequality $xy \leq x \log x - x + e^y$. Using this inequality and the representation $\mathbb{E} \log(1+\zeta) = \mathbb{E}_0\{Z \log(1+\zeta)\}$ with $Z = d\mathbb{P}/d\mathbb{P}_0$, we obtain

$$\mathbb{E} \log(1+\zeta) = \mathbb{E}_0\{Z \log(1+\zeta)\}$$
$$\leq \mathbb{E}_0(Z \log Z - Z) + \mathbb{E}_0(1+\zeta)$$
$$= \mathbb{E}_0(Z \log Z) + \mathbb{E}_0 \zeta - \mathbb{E}_0 Z + 1.$$

It remains to note that $\mathbb{E}_0 Z = 1$ and $\mathbb{E}_0(Z \log Z) = \mathbb{E} \log Z$. $\square$

We now apply this lemma with $\zeta = \varrho(\tilde{\theta}, \theta)/\mathcal{R}(\tilde{\theta}, \theta)$ and show that $\mathbb{E}_0 \zeta = \mathbb{E}_\theta \varrho(\tilde{\theta}, \theta)/\mathcal{R}(\tilde{\theta}, \theta) = 1$. This yields

$$\mathbb{E}_\theta(Z_{I,\theta} \log Z_{I,\theta}) = \mathbb{E} \log Z_{I,\theta} = \mathbb{E} \sum_{t \in I} \log \frac{p(Y_t, f(t))}{p(Y_t, \theta)}$$
$$= \mathbb{E} \sum_{t \in I} \mathbb{E}\left\{\log \frac{p(Y_t, f(t))}{p(Y_t, \theta)} \Big| \mathcal{F}_{t-1}\right\} = \mathbb{E} \Delta_I(\theta)$$

and the result follows. $\square$

This result implies that the bound for the risk of estimation $\mathbb{E}|L_I(\tilde{\theta}_I, \theta)|^r \equiv |N_I \mathbb{E} \mathcal{K}(\tilde{\theta}_I, \theta)|^r$ under the parametric hypothesis can be extended to the nonparametric situation provided that the value $\Delta_I(\theta)$ is sufficiently small. For $r > 0$, define $\varrho(\tilde{\theta}_I, \theta) = |N_I \mathcal{K}(\tilde{\theta}_I, \theta)|^r$. By Theorem 2.4, $\mathcal{R}(\tilde{\theta}_I, \theta) = \mathbb{E}_\theta \varrho(\tilde{\theta}_I, \theta) \leq \mathfrak{r}_r$.

COROLLARY 2.7. *Let (2.5) hold for some $\theta$. For any $r > 0$,*

$$\mathbb{E} \log(1 + |N_I \mathcal{K}(\tilde{\theta}_I, \theta)|^r / \mathfrak{r}_r) \leq \Delta + 1.$$

This result means that in the nonparametric situation under the condition (2.5) with some fixed $\Delta$, the losses $|N_I \mathcal{K}(\tilde{\theta}_I, \theta)|^r$ are stochastically bounded. Note that this result applies even if $\Delta$ is large, however, the bound is proportional to $e^{\Delta+1}$ and grows exponentially with $\Delta$.

2.3. *"Small modeling bias" condition and rate of estimation.* This section briefly comments on relations between the results of Section 2.2 and the usual rate results under smoothness conditions on the function $f(\cdot)$.

Let $n$ be the parameter meaning length of the largest considered historical interval. More precisely, we assume that the function $f(\cdot)$ is smooth in the sense that, for $\theta^* = f(t^\diamond)$ and any $t \geq t^\diamond - n$,

(2.6) $$\mathcal{K}^{1/2}(f(t), \theta^*) \leq (t^\diamond - t)/n.$$

In view of the inequality $\mathcal{K}(\theta, \theta') \asymp |\theta/\theta' - 1|^2$, this condition is equivalent to the usual Lipschitz property of the rescaled function $f(t/n)$. This condition bounds the *bias* of approximating the underlying function $f(t)$ by a constant $f(t^\diamond)$ by $(t^\diamond - t)/n$. The *variance* of the estimate $\tilde{\theta}_I$, for $I = [t, t^\diamond[$, is proportional to $1/(t^\diamond - t)$. The usual "bias-variance trade-off" means the relation "bias$^2 \asymp$ variance," leading to $(t^\diamond - t)^3 \asymp n^2$.

Now note that (2.5) and (2.6) implies

$$\Delta_I(\theta^*) \leq N_I^3/n^2.$$

Therefore, the "small modeling bias" condition $\Delta_I(\theta) \leq \Delta$ is essentially equivalent to "bias-variance trade-off." Moreover, combined with the result of Corollary 2.7, this condition leads to the following classical rate results.

THEOREM 2.8. *Assume (2.6). Select $I$ such that $N_I = cn^{2/3}$, for some $c > 0$. Then, (2.5) holds with $\Delta = c^3$ and, for any $r > 0$,*

$$\log(1 + |N_I \mathcal{K}(\tilde{\theta}_I, \theta)|^r/\mathfrak{r}_r) \leq c^3 + 1.$$

This corresponds to the classical accuracy of nonparametric estimation for the Lipschitz functions [cf. Fan, Farmen and Gijbels (1998)].

**3. Adaptive volatility estimation.** The assumption of time homogeneity is too restrictive in practical applications and it does not allow to fit real data well. In this paper we consider an approach based on the *local parametric assumption*, which means that, for every time moment $t^\diamond$, there exists a historic time interval $[t^\diamond - m, t^\diamond[$ in which the volatility process $\sigma_t$ is nearly constant. Under such a modeling, the main intention is both to describe the interval of homogeneity and to estimate the corresponding value $\sigma_{t^\diamond}$.

Our approach is based on the adaptive choice of the interval of homogeneity for the fixed end point $t^\diamond$. This choice is made by the *local (multiscale) change point detection* (LCP) algorithm described below. The procedure attempts to find this interval from the data by successive testing of the hypothesis of homogeneity. An interval-candidate is accepted if every point is negatively tested on a possible location of a change point. A change point test at a location $\tau < t^\diamond$ compares two different estimates of the parameter $\theta$; one is computed from the most recent interval $[\tau, t^\diamond]$ while the other one is obtained by using another interval $[t', \tau]$ before the possible jump. The procedure is multiscale in the sense that the choice of the other interval $[t', \tau]$ and the critical value of the test depends on the distance of the testing point $\tau$ from the reference point $t^\diamond$. More precisely, let a growing sequence of numbers $N_1 < N_2 < \cdots < N_K$ be fixed. Each $N_k$ means the *scale* parameter describing the length of the historical time interval screened at the step $k$. Define a family $\{\mathcal{I}_k, k = 1, \ldots, K\}$ of nested intervals of the form $\mathcal{I}_k = [t^\diamond - N_k, t^\diamond[$ with the right edge at $t^\diamond$. The procedure starts from the smallest interval $\mathcal{I}_1$ by testing the hypothesis of homogeneity within $\mathcal{I}_1$ against a change point alternative. If the hypothesis is not rejected then we take the next larger interval $\mathcal{I}_k$ and test for a change point. We continue this way until we detect a change point or the largest considered interval $\mathcal{I}_K$ is accepted. If the hypothesis of homogeneity within some interval $\mathcal{I}_k$ is rejected and a change point is detected at a point $\hat{\tau} \in \mathcal{I}_k$, then the estimated interval of homogeneity is defined as the latest accepted interval, that is, $\hat{\mathcal{I}} = \mathcal{I}_{k-1} = [t^\diamond - N_{k-1}, t^\diamond[$, otherwise we take $\hat{\mathcal{I}} = \mathcal{I}_K$. Finally, we define the estimate $\hat{f}(t^\diamond) = \hat{\theta}$ of the volatility parameter $f(t^\diamond) = \sigma_{t^\diamond}^2$ as $\hat{f}(t^\diamond) = \tilde{\theta}_{\hat{\mathcal{I}}}$. The main ingredient of this procedure is the homogeneity test which is described in the next section.

3.1. *Test of homogeneity against a change point alternative.* Let $J$ be a *tested* interval which has to be checked on a possible change point. For carrying out the test, we also assume a larger *testing* interval $I = [t', t''[$ to be fixed. The hypothesis of homogeneity for $J$ means that the observations $R_t$ follow the parametric model with the parameter $\theta$ for $J$ itself and for the larger interval $I$. This hypothesis leads to the parametric log-likelihood $L_I(\theta)$ for the observations $R_t \in I$. We want to test this hypothesis against a change point alternative that the parameter $\theta$ spontaneously changes in some internal point $\tau$ of the interval $J$. Every point $\tau \in J$ splits the interval $I = [t', t''[$ onto two subintervals, $I'' = I''_\tau = [\tau, t''[$ and $I' = I'_\tau = I \setminus I'' = [t', \tau[$ (see Figure 1). The change point alternative means that $f(t) = \theta''$ for $t \in I''$ and $f(t) = \theta'$, for $t \in I''$ for some $\theta'' \neq \theta'$. This corresponds to the log-likelihood $L_{I''}(\theta'') + L_{I'}(\theta')$. The likelihood ratio test statistic for the change point alternative with the change point location at the point $\tau$ is of the form

$$
\begin{aligned}
T_{I,\tau} &= \max_{\theta'', \theta'} \{L_{I''}(\theta'') + L_{I'}(\theta')\} - \max_\theta L_I(\theta) \\
&= \min_\theta \max_{\theta'', \theta'} \{L_{I''}(\theta'', \theta) + L_{I'}(\theta', \theta)\}.
\end{aligned}
\tag{3.1}
$$

For the considered volatility model, this test statistic can be represented in a simple form given by the next lemma.

LEMMA 3.1. *It holds for any interval $I$ and point $\tau \in I$,*

$$T_{I,\tau} = N_{I''} \mathcal{K}(\tilde{\theta}_{I''}, \tilde{\theta}_I) + N_{I'} \mathcal{K}(\tilde{\theta}_{I'}, \tilde{\theta}_I) \tag{3.2}$$

*with $I' = [t', \tau]$ and $I'' = [\tau, t'']$.*

PROOF. By (2.3), minimization in (3.1) with regard to $\theta$ leads to the choice $\theta = \tilde{\theta}_I$. Similarly, maximization with regard to $\theta'$ and $\theta''$ leads to $\theta' = \tilde{\theta}_{I'}$ and $\theta'' = \tilde{\theta}_{I''}$, and the assertion follows. $\square$

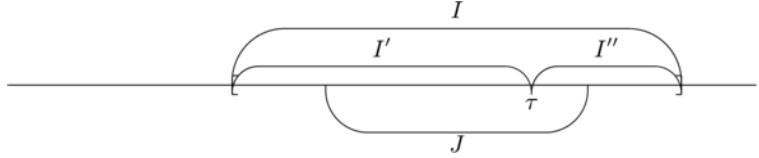

Fig. 1. *Intervals involved in the change point test.*

The change point test for the interval $J$ is defined as the maximum of the test statistics $T_{I,\tau}$ over $\tau \in J$:

$$T_{I,J} = \max_{\tau \in J} T_{I,\tau}. \tag{3.3}$$

The change point test compares this statistic with the critical value $\mathfrak{z}$ which may depend on the intervals $I, J$. The hypothesis of homogeneity is rejected if $T_{I,J} \geq \mathfrak{z}$, in this case the estimate of the change point location is defined as $\hat{\tau} = \arg\max_{\tau \in J_I} T_{I,\tau}$.

REMARK 3.1. The change point alternative suggested above is only one possibility to test the homogeneity assumptions. One can apply many other tests, for example, omnibus tests against polynomials or trigonometric functions [see, e.g., Hart (1998)]. Our choice is motivated by several reasons. First of all, it is simple to implement and does not require a special model estimation under alternative because the alternative reduces to the null hypothesis for two smaller intervals. Secondly, it has a natural interpretation and delivers additional information about the location of the change and the length of the interval of homogeneity. Finally, it was shown in Ingster (1986) [see also Horowitz and Spokoiny (2001)] that a test based on the local constant alternative is powerful against smooth alternatives as well.

3.2. *The multiscale procedure.* This section describes the LCP procedure. The procedure is sequential and consists of $K$ steps corresponding to the given growing sequence of numbers $N_0 < N_1 < \cdots < N_K$. This sequence determines the sequence of nested intervals $\mathcal{I}_0 \subset \mathcal{I}_1 \subset \cdots \subset \mathcal{I}_K$ with the right edge at the point of estimation $t^\diamond : \mathcal{I}_k = [t_k^\diamond, t^\diamond[ = [t^\diamond - N_k, t^\diamond[$ (see Figure 2). This set of intervals leads to the set of estimates $\tilde{\theta}_{\mathcal{I}_k}$, $k = 0, 1, \ldots, K$. Obviously, $N_{\mathcal{I}_k} = N_k$. For conciseness of notation, we also write $\tilde{\theta}_k$ in place of $\tilde{\theta}_{\mathcal{I}_k}$.

The proposed adaptive method chooses an index $\kappa$ of, equivalently, the estimate $\tilde{\theta}_\kappa$ from this set. The procedure is sequential and it successively checked the intervals $\mathcal{I}_0, \mathcal{I}_1, \ldots, \mathcal{I}_k$, on change points.

The interval $\mathcal{I}_0$ is always accepted and the procedure starts with $k = 1$. At every step $k$, every point of the interval $\mathcal{J}_k = \mathcal{I}_k \setminus \mathcal{I}_{k-1}$ is tested as a potential change point due to the procedure from Section 3.1. The testing interval $I = \mathcal{I}_{k+1}$ is applied. $\mathcal{I}_k$ is accepted if the previous interval $\mathcal{I}_{k-1}$ was accepted and the test statistic $T_k \stackrel{\text{def}}{=} T_{\mathcal{I}_{k+1}, \mathcal{J}_k}$ defined by (3.3) does not exceed the critical value $\mathfrak{z}_k$. The latter means that there is no change point detected within $\mathcal{J}_k$. Equivalently, $\mathcal{I}_k$ is accepted if every point is negatively tested on a change point location. The event $\{\mathcal{I}_k$ is rejected$\}$ means that $T_l > \mathfrak{z}_l$ for some $l \leq k$, and hence, a change point has been detected in the first $k$ steps of the procedure at some point within $\mathcal{I}_k$. For every $k$, we define an index $\kappa_k$ corresponding to the largest accepted interval after the first $k$ steps, and $\hat{\theta}_k = \tilde{\theta}_{\kappa_k}$ is the corresponding estimate. The estimate $\tilde{\theta}_k$ and $\hat{\theta}_k$ coincide if no change point is detected at the first $k$ steps. The final estimate is defined as $\hat{\theta} = \hat{\theta}_K$ and corresponds to the largest found interval of homogeneity. The formal definition reads as follows:

$$\kappa = \max\{k \leq K : T_l \leq \mathfrak{z}_l, l = 1, \ldots, k\}, \qquad \hat{\theta} = \tilde{\theta}_\kappa.$$

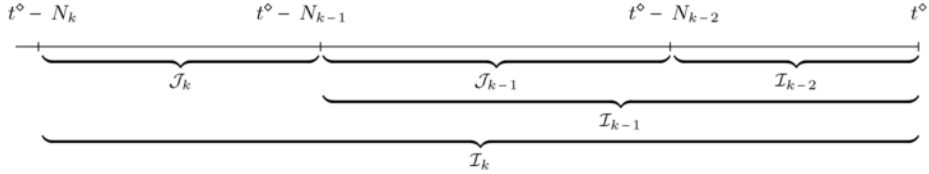

Fig. 2. *The intervals $\mathcal{I}_k$ and $\mathcal{J}_k$ for the LCP procedure.*

The way of choosing the critical value as well as the other parameters of the procedure, like the intervals $\mathcal{I}_k$, is discussed in the next section.

3.3. *Choice of the parameters $\mathfrak{z}_k$ using "propagation" condition.* The "critical value" $\mathfrak{z}_k$ defines the level of significance for the test statistics $T_k = T_{\mathcal{I}_k, \mathcal{J}_k}$. A proper choice of the parameters $\mathfrak{z}_k$ is crucial for the performance of the procedure. We propose in this section one general approach for selecting the $\mathfrak{z}_k$'s, which is similar to the bootstrap idea in the hypothesis testing problem. Indeed, the proposed procedure can be viewed as a multiple test with the scale dependent critical values. We select these values to provide the prescribed performance of the procedure in the parametric situation (under the null hypothesis). In the classical testing approach, the performance of the method is measured by the errors of the first and second kind, and the critical value is selected to ensure the prescribed test level which is the probability of rejecting the null under the null hypothesis. In the considered framework, the null hypothesis means a time-homogeneous model with a constant volatility $\theta^*$. We apply a slightly modified condition on the first kind error which suits better the considered estimation problem. Our primary goal is to select one estimate out of family $\tilde{\theta}_k$, rather than testing on a change point. In the homogeneous situation, our optimal choice, corresponds to the largest interval $\mathcal{I}_K$ leading to the estimate $\tilde{\theta}_K$ with the smallest variability in the considered family (see Theorem 2.4). A "false alarm" means that a nonexisting change point is detected, which leads to selecting an estimate with a larger variability than that of $\tilde{\theta}_K$. Our condition accounts not only for the frequency of false alarms but also at which $k$ step a "false alarm" occurs. Before giving a precise formulation, we mention one important and helpful property of the volatility parametric model $f(\cdot) \equiv \theta^*$: for any intervals $J \subset I$, the distribution of the test statistic $T_{I,J}$ does not depend on the parameter value $\theta^*$. This is a simple corollary of the fact that volatility is a scale parameter of the corresponding parametric family. However, in view of its importance for our study we state it in a separate lemma.

LEMMA 3.2. *Let the return $R_t$ follow the parametric model with the constant volatility parameter $\theta^*$, that is, $R_t^2 = \theta^* \varepsilon_t^2$. Then, for any $J \subset I$ and any $\tau \in J$, the distribution of the test statistics $T_{I,\tau}$ and $T_{I,J}$ under $\mathbb{P}_{\theta^*}$ is the same for all $\theta^* > 0$.*

PROOF. It suffices to notice that for every interval $I$ the estimate $\tilde{\theta}_I$ can be represented under $\mathbb{P}_{\theta^*}$ as
$$\tilde{\theta}_I = N_I^{-1} \sum_{t \in I} Y_t^2 = \theta^* N_I^{-1} \sum_{t \in I} \varepsilon_t^2$$
and for each two intervals $I, I'$, the Kullback–Leibler divergence $\mathcal{K}(\tilde{\theta}_I, \tilde{\theta}_{I'})$ is a function of the ratio $\tilde{\theta}_I / \tilde{\theta}_{I'}$. □

The result of Lemma 3.2 allows us to reduce the parametric null situation to the case of a simple null consisting of one point $\theta^*$, for example, $\theta^* = 1$. The corresponding distribution of the observation under this measure will be denoted by $\mathbb{P}_{\theta^*}$.

For every step $k$, we require that in the parametric situation $f(\cdot) \equiv \theta^*$ the estimate $\hat{\theta}_k$ is sufficiently close to the "oracle" estimate $\tilde{\theta}_k$ in the sense that

$$(3.4) \qquad \mathbb{E}_{\theta^*} |N_k \mathcal{K}(\tilde{\theta}_k, \hat{\theta}_k)|^r \leq \alpha \mathfrak{r}_r$$

for all $k = 1, \ldots, K$ with $\mathfrak{r}_r$, is the parametric risk bound from Theorem 2.4: $\mathbb{E}_{\theta^*} |N_k \mathcal{K}(\tilde{\theta}_k, \theta^*)|^r \leq \mathfrak{r}_r$.

Note that the $\hat{\theta}_k$ differs from $\tilde{\theta}_k$ only if a change point is detected at the first $k$ steps. The usual condition to any change point detector is that such "false alarms" occur with a small probability. Our condition (3.4) has the same flavor but it is a bit stronger. Namely, a false alarm at an early stage of the procedure is more crucial because it results in selecting an estimate with a high variability. Our condition penalizes not only for occurrence of a false alarm but also for the deviation of the selected estimate $\hat{\theta}_k$ from the optimal estimate $\tilde{\theta}_k$. The choice of penalization is motivated by Theorem 2.4. A small value of $N_k \mathcal{K}(\tilde{\theta}_k, \hat{\theta}_k)$ means that $\hat{\theta}_k$ belongs to the confidence set based on the estimate $\tilde{\theta}_k$, that is, $\hat{\theta}_k$ does not differ significantly from $\tilde{\theta}_k$. On the contrary, big values of $N_k \mathcal{K}(\tilde{\theta}_k, \hat{\theta}_k)$ indicate that $\hat{\theta}_k$ differs significantly from $\tilde{\theta}_k$. The choice of the power loss $r$ in the condition (3.4) close to zero leads back to counting the numbers of false alarms. Larger values of $r$ result in the criteria which also accounts for the deviation of $\hat{\theta}_k$ from $\tilde{\theta}_k$. We refer to (3.4) as a "propagation" condition because it ensures that, under homogeneity at every step, the current accepted interval $\hat{I}_k$ extends to $I_k$ with a high probability.

The values $\alpha$ and $r$ in (3.4) are two global parameters. The role of $\alpha$ is similar to the level of the test in the hypothesis testing problem, while $r$ describes the power of the loss function. A specific choice is subjective and depends on the particular application at hand. Taking a large $r$ and small $\alpha$ would result in an increase of the critical values and, therefore, improves the performance of the method in the parametric situation at cost of some loss of sensitivity to parameter changes. Theorem 4.1 presents some upper bounds for the critical values $\mathfrak{z}_k$ as functions of $\alpha$ and $r$ in the form $a_0 \log K + 2\log(N_k/\alpha) + 2r \log(N_K/N_k)$, with some coefficient $a_0$. We see that these bounds linearly depend on $r$ and on $\log \alpha^{-1}$. For our applications to volatility estimation, we apply a relatively small value $r = 1/2$ which makes the procedure more stable and robust against outliers. We also apply $\alpha = 0.2$, although the other values in the range $[0.1, 1]$ lead to very similar results.

The set of conditions (3.4) do not directly define the critical values $\mathfrak{z}_k$. We present below one constructive method for selecting $\mathfrak{z}_k$ to provide the "propagation" conditions (3.4).

3.3.1. *A sequential choice.* Here we present a proposal for a sequential choice of the $\mathfrak{z}_k$'s. Consider the situation after the first $k$ steps of the algorithm. We distinguish between two cases. In the first, change point is detected at some step $l \leq k$, and in the other case no change point is detected. In the first case, we denote by $\mathcal{B}_l$ the event meaning the rejection at the step $l$, that is,

$$\mathcal{B}_l = \{T_1 \leq \mathfrak{z}_1, \ldots, T_{l-1} \leq \mathfrak{z}_{l-1}, T_l > \mathfrak{z}_l\}$$

and $\hat{\theta}_k = \tilde{\theta}_{l-1}$ on $\mathcal{B}_l$, $l = 1, \ldots, k$. The sequential choice of the critical values $\mathfrak{z}_k$ is based on the decomposition

$$(3.5) \qquad |\mathcal{K}(\tilde{\theta}_k, \hat{\theta}_k)|^r = \sum_{l=1}^{k} |\mathcal{K}(\tilde{\theta}_k, \tilde{\theta}_{l-1})|^r \mathbf{1}(\mathcal{B}_l)$$

for every $k \leq K$. Now, we show that the event $\mathcal{B}_l$ only depends on $\mathfrak{z}_1, \ldots, \mathfrak{z}_l$. In particular, the event $\mathcal{B}_1$ means that $T_1 > \mathfrak{z}_1$ and $\hat{\theta}_j = \tilde{\theta}_0$, for all $j \geq 1$. We select $\mathfrak{z}_1$ as the minimal value providing that

$$(3.6) \qquad \max_{k=1,\ldots,K} \mathbb{E}_{\theta^*} |N_k \mathcal{K}(\tilde{\theta}_k, \tilde{\theta}_0)|^r \mathbf{1}(T_1 > \mathfrak{z}_1) \leq \alpha \mathfrak{r}_r / K.$$

Similarly, for every $l \geq 2$, select $\mathfrak{z}_l$ by considering the event $\mathcal{B}_l = \{\kappa = l\}$, meaning that the first false alarm occurs at the step $l$ and $\hat{\theta}_k = \tilde{\theta}_{l-1}$ for all $k > l$. If $\mathfrak{z}_1, \ldots, \mathfrak{z}_{l-1}$ have already been fixed, the event $\mathcal{B}_l$ is only controlled by $\mathfrak{z}_l$ leading to the following condition on $\mathfrak{z}_l$: this is the minimal value that ensures

$$(3.7) \qquad \max_{k \geq l} \mathbb{E}_{\theta^*} |N_k \mathcal{K}(\tilde{\theta}_k, \tilde{\theta}_{l-1})|^r \mathbf{1}(\mathcal{B}_l) = \alpha \mathfrak{r}_r / K.$$

Such a value $\mathfrak{z}_l$ can be found numerically by the Monte Carlo simulations from the parametric model $\mathbb{P}_{\theta^*}$ for any fixed $\theta^*$ [see Lemma 3.2]. It is straightforward to check that such defined $\mathfrak{z}_k$ fulfill (3.4) in view of the decomposition (3.5).

3.3.2. *Examples of choosing the intervals $\mathcal{I}_k$.* To start the procedure, one has to specify the set of intervals $\mathcal{I}_0, \mathcal{I}_1, \ldots, \mathcal{I}_K$. Note, however, that this choice is not a part of the LCP procedure. The method applies to whatever intervals $\mathcal{I}_k$ are selected under condition (MD) (see Section 4). This section presents one example which is at the same time the default choice for our simulation study and applications.

The set $N_0, N_1, \ldots, N_K$ is defined geometrically by the rule $N_k = [N_0 a^k]$ for some fixed value $N_0$ and the growth rate $a > 1$. Such a proposal is motivated by the condition (MD) from the next section. Note also that the sets $\mathcal{J}_k$ do not intersect for different $k$, and every point $\tau \in [t^\diamond - N_k, t^\diamond - N_0]$ is tested as a possible location of the change point at some of the first $k$ steps of the procedure.

Our numerical results (not reported here) indicate that the procedure is quite stable with regard to the choice of the parameters like $N_0$ and $a$. We apply $a = 1.25$. The other values of $a$ in the range 1.1 to 1.3 lead to very similar results. We also apply $N_0 = 5$, which is motivated by our applications to risk management in financial engineering.

**4. Theoretical properties.** This section discusses some useful theoretical properties of the adaptively selected interval of homogeneity $\hat{\mathcal{I}}$ and then of the adaptive volatility estimate $\hat{\theta}$ that corresponds to the selected interval $\hat{\mathcal{I}}$, that is, $\hat{\theta} = \tilde{\theta}_{\hat{\mathcal{I}}}$. Our main "oracle" result claims that the final estimate $\hat{\theta}$ delivers essentially the same quality of estimation as the estimate with the optimal ("ideal") choice of the interval $\mathcal{I}_{k^*}$. It is worth noting that the oracle result does not assume any particular structure of the volatility function $f(t)$. It can be an arbitrary, predictable positive random process. Particular cases include piecewise constant, smooth transition and other models. As shown in Section 2.2, the oracle result automatically ensures the optimal estimation rate under usual smoothness conditions on the function $f(\cdot)$.

The "oracle" result is in its turn, a corollary of two important properties of the procedure: "propagation" under homogeneity and "stability." The first one means that in the nearly homogeneous situation, the procedure would not terminate (no false alarm) with a high probability. In other words, if the parametric (constant) approximation well applies in the interval $\mathcal{I}_k$, then this interval will be accepted with a high probability. The "stability" property ensures that the estimation quality will not essentially deteriorate in the steps "after propagation" when the local constant approximation is not sufficiently accurate. Typically, the procedure terminates in such situations.

The results require some regularity conditions on the growth of the intervals $\mathcal{I}_k$. Namely, we require that the length $N_k$ of $\mathcal{I}_k$ grows exponentially with $k$.

(MD) For some constants $\mathfrak{u}_0, \mathfrak{u}$ with $0 < \mathfrak{u}_0 \leq \mathfrak{u} < 1$, it holds

$$\mathfrak{u}_0 \leq N_{k-1}/N_k \leq \mathfrak{u}.$$

In addition, we assume that the parameter set $\Theta$ satisfies the condition
($\Theta$) for some constant $\mathfrak{a}$ with $0 < \mathfrak{a} \leq 1$, and, for any $\theta_0, \theta \in \Theta$,

$$\mathfrak{a}^2 \leq \theta_0/\theta \leq \mathfrak{a}^{-2}.$$

We start by discussing the behavior of the procedure in the time-homogeneous situation with the constant volatility parameter $\theta^*$. In this case the properties of the resulting estimate $\hat{\theta}$ are guaranteed by the condition (3.4). Our first result claims a possibility of selecting the critical values $\mathfrak{z}_k$ to provide (3.4) and states some upper bounds for the $\mathfrak{z}_k$'s. Similar results can be stated in the local parametric situation when the homogeneity condition $f(t) = \theta^*$ is only fulfilled for some time interval $I$.

4.1. *Behavior under (local) homogeneity.* First, we consider the homogeneous situation with the constant parameter value $f(x) = \theta^*$. Our first result presents an upper bound for the parameters $\mathfrak{z}_k$ which ensures condition (3.4).

THEOREM 4.1. *Assume (MD). Let $f(t) = \theta^*$, for all $t \in \mathcal{I}_K$. Then, there is a constant $a_0 > 0$ depending on $r$ and $\mathfrak{u}_0$, $\mathfrak{u}$ such that the choice*

(4.1) $$\mathfrak{z}_k = a_0 \log K + 2\log(N_k/\alpha) + 2r \log(N_K/N_k)$$

*ensures (3.4), for all $k \leq K$.*

REMARK 4.1. The present result only describes an upper bound for the critical values which will be used for our theoretical study. This upper bound is not used for computing the values $\mathfrak{z}_k$ in practical applications. However, it qualitatively describes how every critical value $\mathfrak{z}_k$ depends on the index $k$ and on the parameter $r, \alpha$.

PROOF OF THEOREM 4.1. Before proving the result of the theorem, we present two useful technical lemmas. The first one shows that the maximum test statistic $T_{I,J}$ is stochastically bounded in a rather strong sense.

LEMMA 4.2. *Let $J, I$ be tested and testing intervals, and $T_{I,J}$ be the test statistic from (3.1). For any other interval $\mathcal{I}$ and any $\mathfrak{z} \geq 1$, it holds*

$$\mathbb{E}_{\theta^*}|N_\mathcal{I}\mathcal{K}(\tilde{\theta}_\mathcal{I}, \theta^*)|^r \mathbf{1}(T_{I,J} > \mathfrak{z}) \leq 2N_J C_r \mathfrak{z}^r e^{-\mathfrak{z}/2},$$

*where $C_r$ is the constant from Theorem 2.4.*

PROOF. Every $\tau \in J$ splits the interval $I$ into $I'$ and $I''$. For any interval $\mathcal{I}$, by Theorem 2.4,

$$\mathbb{E}_{\theta^*}|N_\mathcal{I}\mathcal{K}(\tilde{\theta}_\mathcal{I}, \theta^*)|^r \mathbf{1}(T_{I,\tau} > \mathfrak{z})$$
$$\leq \mathbb{E}_{\theta^*}|N_\mathcal{I}\mathcal{K}(\tilde{\theta}_\mathcal{I}, \theta^*)|^r \{\mathbf{1}(N_{I''}\mathcal{K}(\tilde{\theta}_{I''}, \theta^*) > \mathfrak{z}/2) + \mathbf{1}(N_{I'}\mathcal{K}(\tilde{\theta}_{I'}, \theta^*) > \mathfrak{z}/2)\}$$
$$\leq 2C_r \mathfrak{z}^r e^{-\mathfrak{z}/2}.$$

Now, by definition of $T_{I,J}$,

$$\mathbf{1}(T_{I,J} > \mathfrak{z}) \leq \sum_{\tau \in J} \mathbf{1}(T_{I,\tau} > \mathfrak{z})$$

and the assertion follows. $\square$

Below, we also utilize the metric-like property of the Kullback–Leibler divergence $\mathcal{K}(\theta, \theta')$.

LEMMA 4.3 [Polzehl and Spokoiny (2006), Lemma 5.2]. *Under condition ($\Theta$), it holds that, for every sequence $\theta_0, \theta_1, \ldots, \theta_m \in \Theta$,*

$$\mathcal{K}^{1/2}(\theta_1, \theta_2) \leq \mathfrak{a}\{\mathcal{K}^{1/2}(\theta_1, \theta_0) + \mathcal{K}^{1/2}(\theta_2, \theta_0)\},$$
$$\mathcal{K}^{1/2}(\theta_0, \theta_m) \leq \mathfrak{a}\{\mathcal{K}^{1/2}(\theta_0, \theta_1) + \cdots + \mathcal{K}^{1/2}(\theta_{m-1}, \theta_m)\}.$$

With given constants $\mathfrak{z}_k$, define, for $k > 1$, the random sets

$$\mathcal{A}_k = \{T_k \leq \mathfrak{z}_k\}, \qquad \mathcal{A}^{(k)} = \mathcal{A}_1 \cap \cdots \cap \mathcal{A}_k.$$

Obviously, $\hat{\theta}_k = \tilde{\theta}_k$ on $\mathcal{A}^{(k)}$, for all $k \leq K$, and we have to bound the risk of $\hat{\theta}_k$ on the complement $\bar{\mathcal{A}}^{(k)}$ of $\mathcal{A}^k$. Define $\mathcal{B}_l = \mathcal{A}^{(l-1)} \setminus \mathcal{A}^{(l)}$. The event $\mathcal{B}_l$ means the false alarm at step $l$ and hence, $\kappa = l - 1$. We aim to bound the portion of the risk $\mathbb{E}_{\theta^*}|N_k \mathcal{K}(\tilde{\theta}_k, \hat{\theta}_k)|^r = \mathbb{E}_{\theta^*}|N_k \mathcal{K}(\tilde{\theta}_k, \tilde{\theta}_l)|^r$ associated with the false alarm $\mathcal{B}_l$, for every $l < k$. The definition of $\mathcal{B}_l$ implies that $\mathbf{1}(\mathcal{B}_l) \leq \mathbf{1}(T_{\mathcal{I}_{l+1}, \mathcal{J}_l} > \mathfrak{z}_l)$. Lemma 4.3 and the elementary inequality $(a+b)^{2r} \leq 2^{(2r-1)_+}(a^{2r} + b^{2r})$ with any $a, b \geq 0$ imply

$$\mathbb{E}_{\theta^*}\mathcal{K}^r(\tilde{\theta}_k, \tilde{\theta}_l)\mathbf{1}(\mathcal{B}_l) \leq 2^{(2r-1)_+}\{\mathbb{E}_{\theta^*}\mathcal{K}^r(\tilde{\theta}_k, \theta^*) + \mathbb{E}_{\theta^*}\mathcal{K}^r(\tilde{\theta}_l, \theta^*)\}\mathbf{1}(\mathcal{B}_l).$$

The interval $\mathcal{J}_l$ is of length $N_l$, and by Lemma 4.2,

$$\begin{aligned}
&\mathbb{E}_{\theta^*}\mathcal{K}^r(\tilde{\theta}_k, \tilde{\theta}_l)\mathbf{1}(\mathcal{B}_l) \\
&\quad \leq 2^{(2r-1)_+}\mathfrak{a}^{2r}\{\mathbb{E}_{\theta^*}\mathcal{K}^r(\tilde{\theta}_k, \theta^*) + \mathbb{E}_{\theta^*}\mathcal{K}^r(\tilde{\theta}_l, \theta^*)\}\mathbf{1}(T_{\mathcal{I}_{l+1}, \mathcal{J}_l} > \mathfrak{z}_l) \\
&\quad \leq 2^{(2r-1)_+}\mathfrak{a}^{2r}C_r N_l \mathfrak{z}_l^r e^{-\mathfrak{z}_l/2}(N_k^{-r} + N_l^{-r}).
\end{aligned}$$

This and Theorem 2.4 imply, for every $k \leq K$,

$$\begin{aligned}
\mathbb{E}_{\theta^*}|N_k \mathcal{K}(\tilde{\theta}_k, \hat{\theta}_k)|^r &\leq N_k^r \mathbb{E}_{\theta^*} \sum_{l=1}^{k} \mathcal{K}^r(\tilde{\theta}_k, \tilde{\theta}_l)\mathbf{1}(\mathcal{B}_l) \\
&\leq 2^{(2r-1)_+}\mathfrak{a}^{2r}C_r \sum_{l=1}^{k}\left[1 + \left(\frac{N_k}{N_l}\right)^r\right] 4 N_l \mathfrak{z}_l^r e^{-\mathfrak{z}_l/2}.
\end{aligned}$$

It remains to check using the condition (MD) that the choice $\mathfrak{z}_k = a_0 \log K + 2\log \alpha^{-1} + 2r\log(N_K/N_k) + 2\log N_k$, with properly selected $a_0$, provides the required bound $\mathbb{E}_{\theta^*}|N_k \mathcal{K}(\tilde{\theta}_k, \hat{\theta}_k)|^r \leq \alpha \mathfrak{r}_r$. $\square$

4.2. *Behavior under "small modeling bias" condition.* Now, we extend the previous result to the situation when the parametric assumption is not precisely fulfilled but the deviation from the parametric structure within the considered local model is sufficiently small. At the step $k$, the procedure involves the interval $\mathcal{I}_{k+1}$ used for testing a change point within $\mathcal{J}_k$. Therefore, the deviation from the parametric situation can be measured for the step $k$ by

$$\Delta_k(\theta) \stackrel{\text{def}}{=} \Delta_{\mathcal{I}_{k+1}}(\theta) = \sum_{\mathcal{I}_{k+1}} \mathcal{K}(\theta, f(t))$$

[see (2.4)]. By definition, the modeling bias $\Delta_k(\theta)$ increases with $k$. We suppose that there is a number $k^*$, such that $\Delta_k(\theta)$ is small for some $\theta$ and $k = k^*$, and hence also for all $k \leq k^*$. Consider the corresponding estimate $\hat{\theta}_{k^*}$ obtained after the first $k^*$ steps of the algorithm. Theorem 2.5 implies in this situation the following result.

THEOREM 4.4. *Assume* (MD). *Let $\theta$ and $k^*$ be such that $\mathbb{E}\Delta_{k^*}(\theta) \leq \Delta$, for some $\Delta \geq 0$. Then,*

$$\mathbb{E}\log\left(1 + \frac{|N_{k^*}\mathcal{K}(\tilde{\theta}_{k^*}, \hat{\theta}_{k^*})|^r}{\alpha \mathfrak{r}_r}\right) \leq 1 + \Delta,$$

$$\mathbb{E}\log\left(1 + \frac{|N_{k^*}\mathcal{K}(\tilde{\theta}_{k^*}, \theta)|^r}{\mathfrak{r}_r}\right) \leq 1 + \Delta.$$

4.3. *"Stability after propagation" and "oracle" results.* Due to the "propagation" result, the procedure performs well as long as the "small modeling bias" condition $\Delta_k(\theta) \leq \Delta$ is fulfilled. To establish the accuracy result for the final estimate $\hat{\theta}$, we have to check that the aggregated estimate $\hat{\theta}_k$ does not vary much at the steps "after propagation" when the divergence $\Delta_k(\theta)$ from the parametric model becomes large.

THEOREM 4.5. *Suppose* (MD) *and* ($\Theta$). *Let, for some $k \leq K$, the interval $\mathcal{I}_k$ be accepted by the procedure and hence, $\hat{\theta}_k = \tilde{\theta}_k$. Then, it holds that*

(4.2) $$N_k \mathcal{K}(\hat{\theta}_k, \hat{\theta}_{k+1}) \leq \mathfrak{z}_k.$$

*Moreover, under* (MD), *it holds for every $k'$ with $k < k' \leq K$*

(4.3) $$N_k \mathcal{K}(\hat{\theta}_k, \hat{\theta}_{k'}) \leq \mathfrak{a}^2 c_\mathfrak{u}^2 \bar{\mathfrak{z}}_k$$

*with $c_\mathfrak{u} = (\mathfrak{u}^{-1/2} - 1)^{-1}$ and $\bar{\mathfrak{z}}_k = \max_{l \geq k} \mathfrak{z}_l$.*

REMARK 4.2. An interesting feature of this result is that it is fulfilled with probability one; that is, the control of stability "works" not only with a high probability, but it always applies. This property follows directly from the construction of the procedure.

PROOF OF THEOREM 4.5. If $\mathcal{I}_{k+1}$ is rejected, then $\hat{\theta}_{k+1} = \hat{\theta}_k$ and the assertion (4.2) trivially follows. Now, we consider the case when $\mathcal{I}_{k+1}$ is accepted yielding $\hat{\theta}_k = \tilde{\theta}_k$ and $\hat{\theta}_{k+1} = \tilde{\theta}_{k+1}$. The acceptance of $\mathcal{I}_k$ implies, by definition of the procedure. that $T_k \stackrel{\text{def}}{=} T_{\mathcal{I}_{k+1}, \mathcal{J}_k} \leq \mathfrak{z}_k$ and, in particular, $T_{\mathcal{I}_{k+1}, \tau} \leq \mathfrak{z}_k$, with $\tau = t^\diamond - N_k$ being the left edge of $\mathcal{J}_k$. This yields [see (3.2)] that

$$N_k \mathcal{K}(\tilde{\theta}_k, \tilde{\theta}_{k+1}) \leq \mathfrak{z}_k$$

and the assertion (4.2) is proved.

Now, assumption ($\Theta$) and Lemma 4.3 yield

$$\mathcal{K}^{1/2}(\hat{\theta}_k, \hat{\theta}_{k'}) \leq \mathfrak{a} \sum_{j=k}^{k'-1} \mathcal{K}^{1/2}(\hat{\theta}_j, \hat{\theta}_{j+1}) \leq \mathfrak{a} \sum_{j=k}^{k'-1} (\mathfrak{z}_j/N_j)^{1/2}.$$

The use of assumption (MD) leads to the bound

$$\mathcal{K}^{1/2}(\hat{\theta}_k, \hat{\theta}_{k'}) \leq \mathfrak{a}(\bar{\mathfrak{z}}_k/N_k)^{1/2} \sum_{j=k}^{k'-1} \mathfrak{u}^{(j-k)/2} \leq \mathfrak{a}(1-\sqrt{\mathfrak{u}})^{-1}(\bar{\mathfrak{z}}_k/N_k)^{1/2},$$

which proves (4.3). □

Combination of the "propagation" and "stability" statements implies the main result concerning the properties of the adaptive estimate $\hat{\theta}$.

THEOREM 4.6 ["Oracle" property]. *Let $\mathbb{E}\Delta_k(\theta) \leq \Delta$, for some $\theta \in \Theta$ and $k \leq k^*$. Then, $\hat{\theta}$ is close to the "oracle" estimate $\tilde{\theta}_{k^*}$ in the sense that*

$$\mathbb{E}\log\left(1 + \frac{|N_{k^*}\mathcal{K}(\tilde{\theta}_{k^*}, \hat{\theta}_{k^*})|^r}{\alpha\mathfrak{r}_r}\right) \leq 1 + \Delta,$$

$$N_{k^*}\mathcal{K}(\hat{\theta}_{k^*}, \hat{\theta}) \leq \mathfrak{a}^2 c_\mathfrak{u}^2 \bar{\mathfrak{z}}_{k^*}.$$

The result claims the "oracle" property of the estimate $\hat{\theta}$ in the sense that it belongs with a high probability to a confidence set of the oracle estimate $\tilde{\theta}_{k^*}$, and thus, there is no significant difference between $\hat{\theta}$ and $\tilde{\theta}_{k^*}$.

We also present one corollary about the risk of adaptive estimation for $r = 1/2$. An extension to an arbitrary $r > 0$ is straightforward.

THEOREM 4.7. *Assume* (MD) *and $\mathbb{E}\Delta_{k^*}(\theta) \leq \Delta$, for some $k^*$, $\theta$ and $\Delta$. Then,*

$$\mathbb{E}\log\left(1 + \frac{N_{k^*}^{1/2}\mathcal{K}^{1/2}(\hat{\theta}, \theta)}{\mathfrak{a}\mathfrak{r}_{1/2}}\right) \leq \log\left(1 + \frac{c_\mathfrak{u}\sqrt{\bar{\mathfrak{z}}_{k^*}}}{\mathfrak{r}_{1/2}}\right) + \Delta + \alpha + 1,$$

*where $c_\mathfrak{u}$ is the constant from Theorem 4.5.*

PROOF. By Lemma 4.3, similarly to the proof of Theorem 4.5,

$$\mathfrak{a}^{-1}|N_{k^*}\mathcal{K}(\hat{\theta}, \theta)|^{1/2}$$
$$\leq |N_{k^*}\mathcal{K}(\hat{\theta}_{k^*}, \hat{\theta})|^{1/2} + |N_{k^*}\mathcal{K}(\tilde{\theta}_{k^*}, \hat{\theta}_{k^*})|^{1/2} + |N_{k^*}\mathcal{K}(\tilde{\theta}_{k^*}, \theta)|^{1/2}$$
$$\leq c_\mathfrak{u}\sqrt{\bar{\mathfrak{z}}_{k^*}} + |N_{k^*}\mathcal{K}(\tilde{\theta}_{k^*}, \hat{\theta}_{k^*})|^{1/2} + |N_{k^*}\mathcal{K}(\tilde{\theta}_{k^*}, \theta)|^{1/2}.$$

This, the elementary inequality $\log(1 + a + b) \leq \log(1 + a) + \log(1 + b)$ for all $a, b \geq 0$, Lemma 2.6, Theorem 2.4 and (3.4) yield

$$\mathbb{E}\log(1 + (\mathfrak{a}\mathfrak{r}_{1/2})^{-1}N_{k^*}^{1/2}\mathcal{K}^{1/2}(\hat{\theta}, \theta))$$
$$\leq \log\left(1 + \frac{c_\mathfrak{u}\sqrt{\bar{\mathfrak{z}}_{k^*}}}{\mathfrak{r}_{1/2}}\right)$$

$$+ \mathbb{E} \log \left( 1 + \frac{N_{k^*}^{1/2} \mathcal{K}^{1/2}(\tilde{\theta}_{k^*}, \hat{\theta}_{k^*})}{\mathfrak{r}_{1/2}} + \frac{N_{k^*}^{1/2} \mathcal{K}^{1/2}(\tilde{\theta}_{k^*}, \theta)}{\mathfrak{r}_{1/2}} \right)$$

$$\leq \log \left( 1 + \frac{c_{\mathfrak{u}} \sqrt{\bar{\mathfrak{z}}_{k^*}}}{\mathfrak{r}_{1/2}} \right) + \Delta + \alpha + 1$$

as required. □

REMARK 4.3. Recall that by Theorem 4.4, the "oracle" choice $k^*$ leads to the risk bound for the loss $|N_{k^*} \mathcal{K}(\tilde{\theta}_{k^*}, \theta^*)|^{1/2}$ of the corresponding estimate $\tilde{\theta}_{k^*}$. The adaptive choice states a similar bound but for the loss $|N_{k^*} \mathcal{K}(\hat{\theta}, \theta^*)|^{1/2}/\bar{\mathfrak{z}}_{k^*}^{1/2}$. This means that the accuracy of the adaptive estimate $\hat{\theta}$ is worse by factor $\sqrt{\bar{\mathfrak{z}}_{k^*}}$ which can be considered the payment for adaptation. Due to Theorem 4.1, $\bar{\mathfrak{z}}_{k^*}$ is bounded from above by $a_0 \log K + 2 \log(N_{k^*}/\alpha) + 2r \log(N_K/N_{k^*})$. Therefore, the risk of the adaptive estimate corresponds to the best possible risk among the family $\{\tilde{\theta}_k\}$, for the choice $k = k^*$ up to a logarithmic factor in the sample size. Lepski, Mammen and Spokoiny (1997) established a similar result in the regression setup for the pointwise adaptive Lepski procedure. Combining the result of Theorem 4.7 with Theorem 2.8 yields the rate of adaptive estimation $(n^{-1} \log n)^{1/(2+d)}$ under Lipschitz smoothness of the function $f$ and the usual design regularity [see Polzehl and Spokoiny (2006) for more details]. It was shown by Lepski (1990) that in the problem of pointwise adaptive estimation this rate is optimal and cannot be improved by any estimation method. This gives an indirect proof of the optimality of our procedure. The factor $\bar{\mathfrak{z}}_{k^*}$ in the accuracy of estimation cannot be removed or reduced in the rate because otherwise the similar improvement would appear in the rate of estimation.

4.4. *Switching regime model.* A *switching regime* model is described by a sequence $\nu_1 < \nu_2 < \cdots$ of Markov moments with respect to the filtration $(\mathcal{F}_t)$ and by values $\theta_1, \theta_2, \ldots$, where each $\theta_j$ is $\mathcal{F}_{\nu_j}$-measurable. By definition, $\sigma_t^2 = f(t) = \theta_j$, for $\nu_j \leq t < \nu_{j+1}$, and $\sigma_t$ is constant for $t < \nu_1$. This is an important special case of the model (2.1). It is worth mentioning that any volatility process $\sigma_t$ can be approximated by a switching regime model. For this special case, the above procedure has a very natural interpretation. When estimating at the point $t^\diamond$, we search for a largest interval of the form $[t^\diamond - m, t^\diamond[$ which does not contain a change point. More precisely, with a giving sequence of interval-candidates $\mathcal{I}_k = [t^\diamond - N_k, t^\diamond[$, we are looking for the largest homogeneous interval among them. This is done via successive testing for a change point within the intervals $\mathcal{J}_k = [t^\diamond - N_k, t^\diamond - N_{k-1}[$.

The construction of the procedure automatically provides the prescribed risk level associated with the first kind error (a false alarm). In this section, we aim to show that the procedure ensures a near-optimal quality of change point detection. The quality (sensitivity) of a change point procedure is usually measured by the mean delay between the occurrence of the change points and its detection. To study this property of the proposed method, we consider the case of estimation at a point $t^\diamond$ next after a change point $\nu$. The "ideal" choice $\mathcal{I}_{k^*}$ among $\mathcal{I}_1, \ldots, \mathcal{I}_K$ is obviously the largest one which does not contain $\nu$. Theorem 4.7 claims that the procedure accepts with a high probability all the intervals $\mathcal{I}_k$ for which the testing interval $\mathcal{I}_{k+1}$ does not contain the point of change $\nu$. This particularly implies that the quality of estimation of $\theta_{t^\diamond}$ by our adaptive procedure is essentially the same as if we knew the latest change point $\nu$ a priori. Now we additionally show that the procedure rejects with a high probability the first interval $\mathcal{I}_{k^*+1}$ which contains the point of change $\nu$ provided that the magnitude of the change is sufficiently large. This fact can be treated as the sensitivity of the procedure to the changes of regime.

In our study, we assume that the changes occur not too often, and there is exactly one change within $\mathcal{J}_{k^*+1}$ and moreover, within the larger interval $\mathcal{I}_{k^*+2}$ which is used as the testing interval for $\mathcal{J}_{k^*+1}$. Let $\theta'$ be the value of the parameter before the change and $\theta''$ after it. The point $\tau$ splits the interval $I = \mathcal{I}_{k^*+2}$ into two homogeneous intervals. $f(t) = \theta''$ for $t \in I'' = [\tau, t^\diamond[$, while $f(t) = \theta_\ell$ within the complementary interval $t \in I \setminus I''$. Define $c_1 = N_{k^*}/N_{k^*+2}$, $c_2 = N_{k^*+1}/N_{k^*+2}$. By condition (MD), $c_1 \geq \mathfrak{u}_0^2$ and $c_2/c_1 \geq \mathfrak{u}^{-1}$. The length $t^\diamond - \tau$ of the interval $[\tau, t^\diamond[$ fulfills $c_1 \leq (t^\diamond - \tau)/N_{k^*+2} \leq c_2$. Based on these considerations, define the following measure of change from $\theta'$ to $\theta''$:

$$(4.4) \qquad d^2(\theta', \theta'') = \inf_\theta \inf_{c \in [c_1, c_2]} \{(1-c)\mathcal{K}(\theta', \theta) + c\mathcal{K}(\theta'', \theta)\}.$$

The following simple bound can be useful for bounding the distance $d^2(\theta', \theta'')$.

LEMMA 4.8. *There is a constant $\mathfrak{b} > 0$ depending on $c_1$ and $c_2$ only such that*

$$d^2(\theta', \theta'') \geq \mathfrak{b}(\theta'/\theta'' - \theta''/\theta')^2.$$

PROOF. For every fixed $\theta', \theta'', \theta$, the expression $(1-c)\mathcal{K}(\theta', \theta) + c\mathcal{K}(\theta'', \theta)$ is a linear function of $c$. Therefore, its minimum with regard to $c$ is attained in one of the edge points $c_1, c_2$, and it suffices to check the assertion only for these two values of $c$. Now, the assertion follows directly from the definition of the Kullback–Leibler distance $\mathcal{K}(\theta', \theta)$ as a smooth function of the ratio $\theta'/\theta$ with $\mathcal{K}(\theta, \theta) \equiv 0$. □

We aim to show that if the contrast $d(\theta', \theta'')$ is sufficiently large then the test statistic $T_{k^*+1}$ will be large as well, yielding that the interval $\mathcal{I}_{k^*+1}$ will be rejected with a high probability.

THEOREM 4.9. *Let $f(t) = \theta'$ before the change point at $\nu$ and $f(t) = \theta''$ after it. If, for some $\mathfrak{z} > 0$,*

(4.5) $$N_{k^*+2} d^2(\theta', \theta'') \geq 2\mathfrak{a}^2 (\mathfrak{z}_{k^*+1} + \mathfrak{z}) \cdot m,$$

*then*

$$\mathbb{P}(\mathcal{I}_{k^*+1} \text{ is not rejected}) \leq 4e^{-\mathfrak{z}}.$$

PROOF. Let $\nu$ be the location of the change within $\mathcal{J}_{k^*+1} = \mathcal{I}_{k^*+1} \setminus \mathcal{I}_{k^*}$. It suffices to show that, under the conditions of the theorem, the corresponding test statistic $T_{I,\nu}$ exceeds with a high probability the value $\mathfrak{z}_{k^*+1}$. This would ensure that the interval $\mathcal{I}_{k^*+1}$ is rejected. The point $\nu$ splits the testing interval $I = \mathcal{I}_{k^*+2}$ into two subintervals $I'$ and $I''$, and within each of the intervals $I'$ and $I''$ the function $f(t)$ is constant. $f(t) \equiv \theta'$ for $t \in I'$, and $f(t) \equiv \theta''$ for $t \in I''$. Let a value $\mathfrak{z} > 0$ be fixed. Introduce the event

$$\mathcal{A}(\mathfrak{z}) = \mathbf{1}\{N_{I'}\mathcal{K}(\tilde{\theta}_{I'}, \theta') \leq \mathfrak{z},\ N_{I''}\mathcal{K}(\tilde{\theta}_{I''}, \theta'') \leq \mathfrak{z}\}.$$

By Theorem 2.2,

$$\mathbb{P}(\mathcal{A}(\mathfrak{z})) \geq 1 - 4e^{-\mathfrak{z}}.$$

We now consider $\mathfrak{z}$ such that (4.5) holds and show that $T_{I,\nu} > \mathfrak{z}_{k^*+1}$ on $\mathcal{A}(\mathfrak{z})$. By definition, it holds on the set $\mathcal{A}(\mathfrak{z})$ that $N_{I'}\mathcal{K}(\tilde{\theta}_{I'}, \theta') \leq \mathfrak{z}$ and $N_{I''}\mathcal{K}(\tilde{\theta}_{I''}, \theta'') \leq \mathfrak{z}$. By Lemma 4.3,

$$\mathcal{K}^{1/2}(\theta', \tilde{\theta}_I) \leq \mathfrak{a}\mathcal{K}^{1/2}(\tilde{\theta}_{I'}, \theta') + \mathfrak{a}\mathcal{K}^{1/2}(\tilde{\theta}_{I'}, \tilde{\theta}_I)$$
$$\leq \mathfrak{a}(\mathfrak{z}/N_{I'})^{1/2} + \mathfrak{a}\mathcal{K}^{1/2}(\tilde{\theta}_{I'}, \tilde{\theta}_I).$$

Hence,

$$\mathcal{K}(\theta', \tilde{\theta}_I) \leq 2\mathfrak{a}^2 \mathfrak{z}/N_{I'} + 2\mathfrak{a}^2 \mathcal{K}(\tilde{\theta}_{I'}, \tilde{\theta}_I)$$

and

$$\mathcal{K}(\tilde{\theta}_{I'}, \tilde{\theta}_I) \geq (2\mathfrak{a}^2)^{-1} \mathcal{K}^{1/2}(\theta', \tilde{\theta}_I) - \mathfrak{z}/N_{I'}.$$

Similarly,

$$\mathcal{K}(\tilde{\theta}_{I''}, \tilde{\theta}_I) \geq (2\mathfrak{a}^2)^{-1} \mathcal{K}^{1/2}(\theta'', \tilde{\theta}_I) - \mathfrak{z}/N_{I''}.$$

Now, by definition of $T_{I,\nu}$ [see (3.2)]

$$T_{I,\nu} = N_{I'}\mathcal{K}(\tilde{\theta}_{I'}, \tilde{\theta}_I) + N_{I''}\mathcal{K}(\tilde{\theta}_{I''}, \tilde{\theta}_I)$$
$$\geq (2\mathfrak{a}^2)^{-1}\{N_{I'}\mathcal{K}(\theta', \tilde{\theta}_I) + N_{I''}\mathcal{K}(\theta'', \tilde{\theta}_I)\} - \mathfrak{z}$$
$$= (2\mathfrak{a}^2)^{-1} N_{k^*+2}\{c\mathcal{K}(\theta', \tilde{\theta}_I) + (1-c)\mathcal{K}(\theta'', \tilde{\theta}_I)\} - \mathfrak{z}$$

with $c = N_{I'}/N_{k^*+2}$. This and the definition of $d(\theta', \theta'')$ [see (4.4)] yields on $\mathcal{A}(\mathfrak{z})$

$$T_{I,\nu} \geq \frac{N_{k^*+2}}{(2\mathfrak{a}^2)} d^2(\theta', \theta'') - \mathfrak{z}.$$

The theorem assertion follows. $\square$

The result of Theorem 4.9 delivers some additional information about the sensitivity of the proposed procedure to changes in the volatility parameter. One possible question is about the minimal delay $m^*$ between the change point $\nu$ and the first moment $t^\diamond$ when the procedure starts to detect this change. Due to Theorem 4.9, the change will be "detected" with a high probability if (4.5) meets. With fixed $\theta' \neq \theta''$, condition (4.5) is fulfilled if $m^*$ is larger than a prescribed constant, that is, we need only a finite number of observations to detect a change point. In general, $m^*$ should be of order $d^{-2}(\theta', \theta'') \asymp |\theta' - \theta''|^{-2}$, if the size of the change becomes small. All these issues are in agreement with the theory of change point detection [see, e.g., Pollak (1985), Csörgő and Horváth (1997) and Brodskij and Darkhovskij (1993)] and with our numerical results from Section 5.

**5. Simulated results and applications.** This section illustrates the performance of the proposed local change point detection (LCP) procedure by means of some simulated and real data sets. We aim to show that the theoretical properties of the method derived in the previous section are confirmed by the numerical results. We focus especially on the two main features of the method: stability under homogeneity and sensitivity to changes of volatility.

5.1. *Some simulated examples.* Three different jump processes are simulated, whose relative jump magnitude is 3.00, 2.00 and 1.75, respectively. Each process is simulated and estimated 1000 times, and the median and the quartiles of the estimates are plotted in Figure 3. We show the results for the final estimate $\hat{\theta}$ and for the length of the selected interval $\hat{\mathcal{I}}$. One can see that, if the size of the change is large enough, the procedure performs as if the location of the change were known. As one can expect, the sensitivity of the change point detection decreases when the magnitude of the jump becomes smaller. However, the accuracy of the volatility estimate remains rather good even for small jumps that corresponds to our theoretical results.

The algorithm proposed in this paper is compared with the LAVE procedure from Mercurio and Spokoiny (2004) with the optimized tuning parameters $\gamma = 0.5$, $M = 40$ and $\mathfrak{z} = 2.40$. Figure 4 shows the quartiles of estimation for the two approaches for the model with the relative jump magnitude 3. One can see that the new procedure outperforms LAVE both with respect to the variance and to the bias of the estimator, especially for the points immediately after the changes.

Our simulation study has been done for the conditional normal model (2.1). We mentioned in Section 2.1 that this assumption is questionable as far as the real financial data is considered. To gain an impression about the robustness of the method against violation of the normality assumption, we also simulated using i.i.d. innovations from the $t_5$-distribution with five degrees of freedom. The results are shown in Figure 5. As one can expect, they are slightly worse than in the case of normal innovations, however, the procedure continues to work in a quite reasonable way. The sensitivity of the procedure remains as good as with normal innovations, but a probability to reject a homogeneous interval became larger. This results in a higher variability of the estimated volatility.

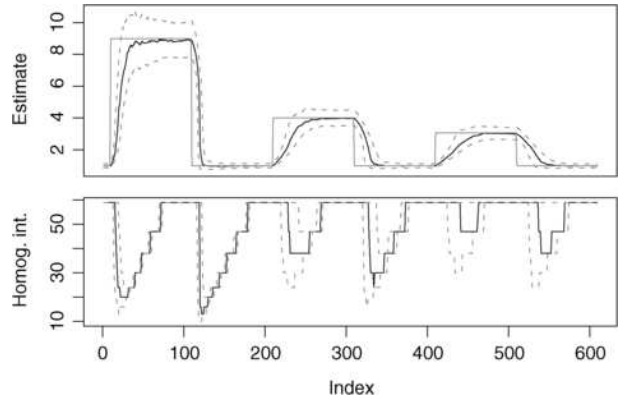

Fig. 3. *A process with Gaussian innovations and jumps of different magnitudes. Top panel: jump process (thin solid line), pointwise median (solid line) and quartiles (dashed lines) for the estimates $\hat{\theta}_t$. Bottom panel: length of the selected interval $\hat{\mathcal{I}}_t$ (solid line) and its quartiles (dashed lines). The results were obtained with parameters $r = 0.5$ and $\alpha = 0.2$ and interval lengths $5, 7, 10, 13, 16, 20, 24, 30, 38, 47, 59, 73, 92$ points.*

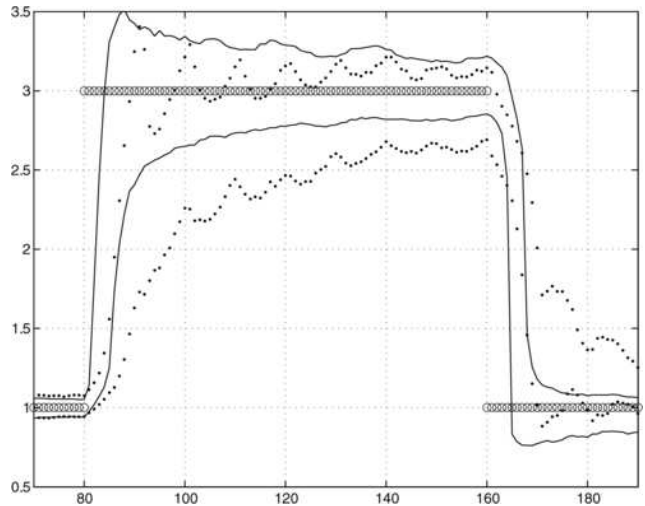

Fig. 4. *Comparison of the proposed estimator with the one from Mercurio and Spokoiny (2004) for change point model with $\theta/\theta' = 3$. Quartiles of $\hat{\theta}$ for the LCP method (solid lines) and for the LAVE method from Mercurio and Spokoiny (2004) (dotted lines) against the true volatility (thick line).*

5.2. *Volatility estimation for different exchange rate data sets.* The volatility estimation is performed on a set of nine exchange rates, which are available from the web page of the US Federal Reserve. The data sets represent daily exchange rates of the US Dollar (USD) against the following currencies: Australian Dollar (AUD), British Pound (GBP) Canadian Dollar (CAD), Danish Krone (DKR), Japanese Yen (JPY), Norwegian Krone (NKR), New Zealand Dollar (NZD), Swedish Krone (SKR) and Swiss Franc (SFR). The period under consideration goes from January 1, 1990 to April 7, 2000. For each time series we have 2583 observations. All selected time series display excess kurtosis and volatility clustering.

Figure 6 show the GBP/USD exchange rate returns together with the volatility estimated with the default parameters. The results of the estimation are in accordance with the data, and the procedure seems to recognize changes in the underlying volatility process quickly.

The assumption of local homogeneity leads to the constant forecast $\hat{\sigma}_t^2$ of the volatility $\sigma_{t+h}$ for a small or moderate time horizon $h$. This results in

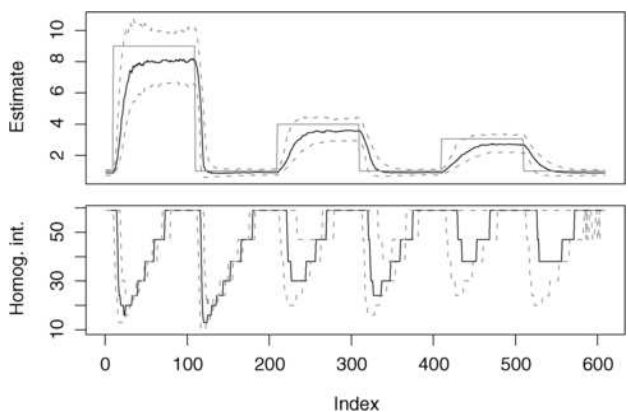

FIG. 5. *Estimation results with respect to jump processes with jumps of different magnitudes. The results are obtained with tuning parameters $r = 0.5$ and $\alpha = 0.2$ and interval lengths $5, 7, 10, 13, 16, 20, 24, 30, 38, 47, 59, 73, 92$ points. The conditional distribution is scaled student $t_5$ with five degrees of freedom.*

the following forecast of the conditional variance of the aggregated returns $R_{t+1} + \cdots + R_{t+h}$:

$$V_{t,h}^{\text{LCP}} := h\hat{\sigma}_t^2.$$

In order to assess the performance of the proposed algorithm, we compare its forecasting ability with that of the $\text{GARCH}(1,1)$ model, which represents one of the most popular parameterizations of the volatility process of financial time series. The $\text{GARCH}(1,1)$ model is described by the following equations:

$$R_t = \sigma_t \varepsilon_t, \qquad \sigma_t^2 = \omega + \alpha R_{t-1}^2 + \beta \sigma_{t-1}^2,$$
$$\alpha > 0, \qquad \beta > 0, \qquad \alpha + \beta < 1, \qquad \varepsilon_t \sim \mathcal{N}(0,1) \qquad \forall t.$$

The $h$-step ahead variance forecast of the $\text{GARCH}(1,1)$ is given by

$$\sigma_{t+h|t}^{2,\text{GARCH}} := \mathbb{E}_t R_{t+h}^2 = \bar{\sigma}^2 + (\alpha + \beta)^h (\sigma_t^2 - \bar{\sigma}^2),$$

where $\bar{\sigma}$ represents the unconditional volatility and $\mathbb{E}_t \xi$ means $\mathbb{E}(\xi|\mathcal{F}_t)$ [see Mikosch and Starica (2000a)]. Since the returns are conditionally uncorrelated, the conditional variance of the aggregated returns is given by the sum of the conditional variances:

$$V_{t,h}^{\text{GARCH}} := \mathbb{E}_t(R_{t+1} + \cdots + R_{t+h})^2 = \sum_{k=1}^{h} \mathbb{E}_t R_{t+k}^2 = \sum_{k=1}^{h} \sigma_{t+h|t}^{2,\text{GARCH}}.$$

The assumption of constant parameters for a $\text{GARCH}(1,1)$ model over a time interval of the considered length of about 2500 time points can be too restrictive. We therefore considered a scrolling estimate, that is, for every date, the preceding 1000 observations are used for estimation of the GARCH parameters, and then the estimated parameters are used to forecast the variance at different horizons. This method is nonadaptive in the choice of the observation window, but it takes advantage of a more flexible GARCH-modeling. The LCP algorithm suggested in this paper applies a very simple local constant modeling but benefits from a data-driven choice of the interval of homogeneity.

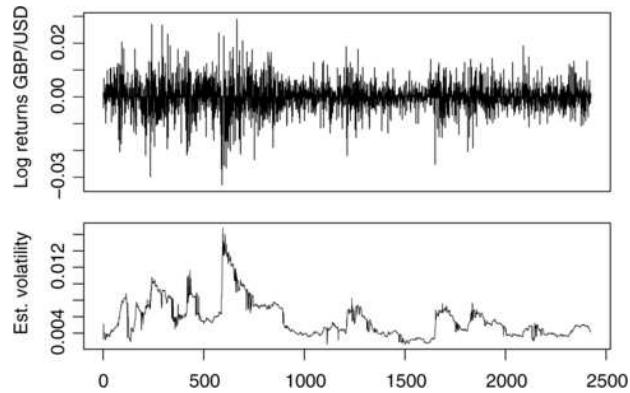

Fig. 6. *Returns and estimated volatility for the GBP/USD exchange rate.*

The quality of forecasting is measured by comparing the forecasts $V_{t,h}^{\text{LCP}}$, respectively, $V_{t,h}^{\text{GARCH}}$ with the realized volatility

$$\bar{V}_{t,h} := R_{t+1}^2 + \cdots + R_{t+h}^2.$$

We apply the following mean square root error criterion (MSqE) for an interval $I$:

$$\text{MSqE}_I = \sum_{t \in I} |V_{t,h}^{\text{LCP}} - \bar{V}_{t,h}|^{1/2} \Big/ \sum_{t \in I} |V_{t,h}^{\text{GARCH}} - \bar{V}_{t,h}|^{1/2}.$$

The MSqE is considered instead of the more common MSE for robustness reasons. Actually, in this way outliers are prevented from having a strong influence on the results. The MSqE is computed for six nonoverlapping intervals of 250 observations, and the results are shown in Table 1. One can observe that both methods are comparable and that the relative performance depends on the particular situation at hand. For periods with stable volatility the LCP forecast is clearly better, but for periods with high volatility variation the GARCH method is slightly preferable.

5.3. *Analysis of standardized returns.* Our model (2.1) assumes the standard normal innovations $\varepsilon_t$. Many empirical researches argued that this assumption is too strong and often violated [see, e.g., McNeil and Frey (2000)]. Here, we briefly discuss this issue by looking at the standardized returns $\hat{\xi}_t = R_t/\hat{\sigma}_t$. The first observation is that even after standardization by the estimated variance, the density of standardized returns $\hat{\xi}_t$ still displays tails which are fatter than the normal. We illustrate this effect in Figure 7 where the kernel estimate of the density of standardized returns $R_t/\hat{\sigma}_t$ is plotted against the normal density and the scaled student $t_5$ density with five degrees of freedom. One can observe that the $t$-distribution delivers a much better approximation to the empirical density of returns.

The volatility clustering effect, though, disappears after standardization and autocorrelations of squared returns are not significant any more (see Figure 8 for the case of GBP/USD returns). The other exchange rate examples deliver similar results. A short conclusion of this empirical study is that the standardized returns can be treated as i.i.d. random variables with a distribution whose tails are fatter than the ones of the normal distribution.

TABLE 1
*Forecasting performance (MSqE) of LCP relative to* GARCH(1, 1)
*on six consecutive time periods of 250 observations each*

| | | | | | | | |
|---|---|---|---|---|---|---|---|
| CAD | $h = 1$ | 0.994 | 0.983 | 0.833 | 0.967 | 1.022 | 0.998 |
| | $h = 5$ | 0.941 | 0.999 | 0.720 | 0.994 | 1.105 | 1.009 |
| | $h = 10$ | 0.862 | 1.038 | 0.645 | 0.960 | 1.149 | 0.999 |
| DKR | $h = 1$ | 0.881 | 0.924 | 0.844 | 0.979 | 0.976 | 1.013 |
| | $h = 5$ | 0.849 | 0.968 | 0.802 | 1.035 | 0.987 | 1.007 |
| | $h = 10$ | 0.870 | 0.971 | 0.691 | 1.053 | 0.986 | 0.989 |
| JPY | $h = 1$ | 0.931 | 0.987 | 0.892 | 1.004 | 1.021 | 0.992 |
| | $h = 5$ | 0.876 | 1.006 | 0.858 | 1.002 | 1.032 | 0.998 |
| | $h = 10$ | 0.889 | 0.978 | 0.828 | 1.033 | 1.061 | 1.001 |
| AUD | $h = 1$ | 0.973 | 0.919 | 0.895 | 1.017 | 1.022 | 0.993 |
| | $h = 5$ | 0.966 | 0.943 | 0.877 | 1.012 | 0.967 | 0.959 |
| | $h = 10$ | 0.932 | 0.958 | 0.887 | 1.032 | 1.023 | 0.990 |
| GBP | $h = 1$ | 0.874 | 0.969 | 0.904 | 1.029 | 0.947 | 0.960 |
| | $h = 5$ | 0.814 | 0.960 | 0.914 | 1.090 | 0.941 | 0.952 |
| | $h = 10$ | 0.775 | 0.890 | 0.884 | 1.087 | 0.972 | 0.949 |
| NZD | $h = 1$ | 0.845 | 0.941 | 0.928 | 1.042 | 0.987 | 0.700 |
| | $h = 5$ | 0.816 | 0.918 | 0.913 | 1.065 | 1.002 | 0.657 |
| | $h = 10$ | 0.742 | 0.984 | 0.884 | 1.095 | 1.013 | 0.632 |

5.4. *Application to value at risk.* The Value at Risk (VaR) measures the extreme loss of a portfolio over a predetermined holding period with a prescribed confidence level $1 - \alpha$. This problem can be reduced to computing the quantiles of the distribution of aggregated returns [see, e.g., Fan and Gu (2003)] for a recent overview of this topic.

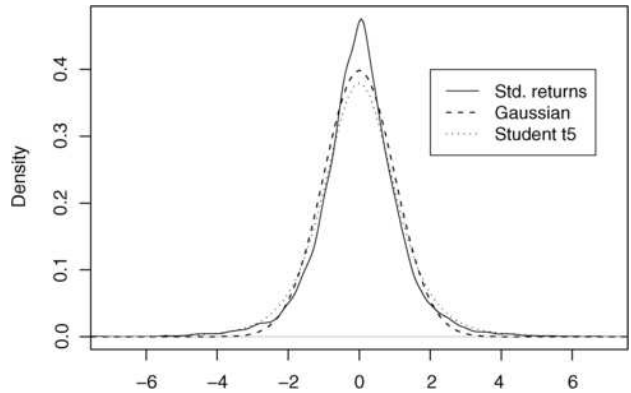

FIG. 7. *Kernel density estimate of exchange rate returns JPY/USD, normal density and scaled student's $t_5$ density with five degrees of freedom.*

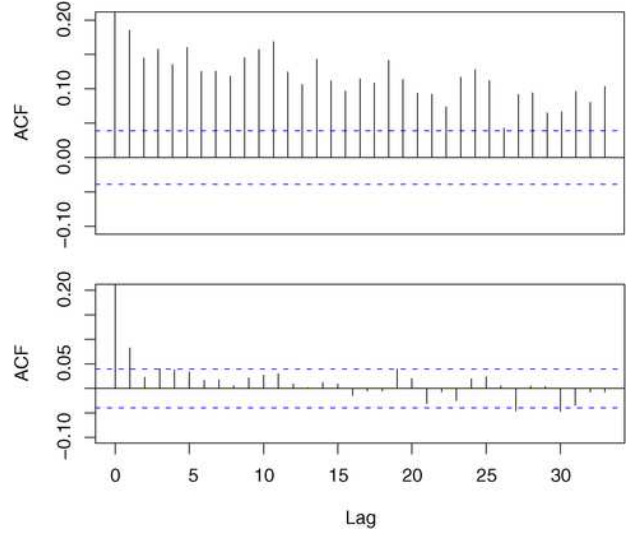

FIG. 8. *ACF of the absolute GBP/USD returns (upper plot) and of the standardized absolute GBP/USD returns (lower plot). Dotted straight line denotes the 95% significance level.*

Our modeling approach can easily be adapted to the VaR problem. Namely, one may forecast the 1% and 5% quantile of the next return $R_{t+1}$ and of the aggregated returns $R_{t+1} + \cdots + R_{t+h} = \log(S_{t+h}/S_t)$, for each date $t$, in the following way. The volatility parameter $\hat{\sigma}_t$ is estimated from the historical data $R_s$, for $s \leq t$, and one can consider different distributions for the innovations $\varepsilon_t$. In our study, we compare the Gaussian, the scaled student $t_5$-distribution with five degrees of freedom and the empirical distribution $\hat{F}_t$ of the past empirical innovations $\hat{\xi}_s$ for $s \leq t$, that is,

$$R_{t+h} = \hat{\sigma}_t \xi_{t+h} \quad \text{with } \xi_{t+h} \sim \mathcal{N}(0,1) \quad \text{or} \quad \sqrt{5/3}\xi_{t+h} \sim t_5 \quad \text{or} \quad \xi_{t+h} \sim \hat{F}_t.$$

Similar approaches have been applied in McNeil and Frey (2000) with the use of the GARCH(1,1) model for estimating the volatility and extreme value theory for evaluating the distribution of returns, while Eberlein and Prause (2002) assume the Generalized Hyperbolic Distribution for the innovations.

In order to better interpret the results, we notice that the scaled $t_5$ distribution has higher 5%-quantiles than the ones of the Gaussian at any of the considered horizons and lower 1%-quantiles. Therefore, the Gaussian distribution of innovations is more conservative for 5%-quantiles, while the opposite is true for 1%-quantiles.

We apply the procedure to the set of nine exchange rates with about 2500 observations in each one. The frequency of overshooting the predicted quantile for the given realizations of the returns is given in Table 2. The first 500 observations in every time series are taken as presample for estimating the parameters. Notice, that for the five and ten day horizon, overlapping intervals of data are used as in Fan and Gu (2003).

According to the requirement of the regulators [BIS (1996)], a bank has to determine its capital requirements in order to cover from market risk proportionally to the 1% quantile of the distribution of the portfolio losses over a ten-day horizon. Internal models calculating this quantile are regularly monitored. The coefficient of proportionality is set to 3 for models whose performance is satisfactory (green zone), and it can be increased up to 4 by a discretionary judgment of the regulators for models which appear to estimate the quantile imprecisely (yellow zone). If the model performance is considered very poor, the coefficient is automatically increased to 4 (red zone).

The official criterion for the evaluation of an internal model is the statistical significance of the 1% quantile estimates of the portfolio loss distribution over a one-day horizon. The prescribed procedure, called backtesting, checks whether the observed frequency of days out of the last 250 for which the losses were larger than the value computed by the prescribed VaR procedure does not significantly deviate from the nominal level 0.01 [see Deutsch (2001)]. Every procedure is classified as green, yellow and red. The green zone means that the empirical frequency is in agreement with the nominal probability 0.01. The yellow zone begins at the point, such that the probability of exceptions for the tested VaR procedure exceeds the value 0.01 with a 95% confidence interval. One can easily verify that such probability corresponds to 5 or more exceptions out of 250 days, that is, the frequency of exceptions equals 2%. Similarly, the red zone corresponds to the 99.99% level, evidence that the tested procedure does not provide the required probability of exceptions. For a sample of 250 observations, this corresponds to 10 exceptions, or equivalently, 4% frequency of overshooting the VaR value.

TABLE 2
Percentage of overshooting the prescribed VaR level for six series of exchange rates for nominal quantile levels 1% and 5%, three different distributions of innovations and time horizon $h = 1, 5, 10$

|  |  | Gaussian | | | Student 5 | | | Empirical | | |
|---|---|---|---|---|---|---|---|---|---|---|
|  | $h$ | 1 | 5 | 10 | 1 | 5 | 10 | 1 | 5 | 10 |
| 1% quantile | AUD | 2.0 | 1.7 | 1.0 | 1.6 | 1.5 | 0.8 | 0.6 | 1.2 | 0.3 |
|  | CAD | 1.9 | 2.2 | 1.9 | 1.4 | 1.8 | 1.6 | 0.6 | 1.2 | 1.0 |
|  | DKR | 1.1 | 1.3 | 1.0 | 0.5 | 0.9 | 0.8 | 0.2 | 0.6 | 0.4 |
|  | GBP | 1.6 | 1.6 | 1.2 | 1.2 | 1.3 | 1.1 | 0.3 | 0.8 | 0.5 |
|  | JPY | 0.8 | 0.6 | 0.5 | 0.5 | 0.4 | 0.4 | 0.1 | 0.1 | 0.1 |
|  | NZD | 2.3 | 1.7 | 1.4 | 1.8 | 1.4 | 1.1 | 0.8 | 0.9 | 0.4 |
| 5% quantile | AUD | 5.4 | 5.5 | 4.9 | 6.3 | 5.6 | 5.1 | 4.6 | 4.3 | 5.8 |
|  | CAD | 5.5 | 6.7 | 7.1 | 6.2 | 6.9 | 7.2 | 4.5 | 5.0 | 7.6 |
|  | DKR | 5.0 | 5.6 | 5.4 | 6.0 | 5.7 | 5.7 | 4.2 | 4.0 | 6.2 |
|  | GBP | 5.1 | 5.7 | 5.5 | 5.9 | 5.8 | 5.7 | 4.3 | 4.1 | 6.1 |
|  | JPY | 4.0 | 3.7 | 4.0 | 5.0 | 3.9 | 4.1 | 3.4 | 2.3 | 4.6 |
|  | NZD | 5.0 | 5.5 | 5.4 | 5.4 | 5.7 | 5.6 | 4.1 | 4.3 | 6.1 |

The comparison of these requirements with our results presented in Table 2 shows that, on average, none of the procedures we tried are in the red zone, and that the procedure using empirical distribution function for the residuals is always in the green zone. The use of the student $t_5$ distribution also allows us to get the green zone results for most of the examples, while the procedure with Gaussian innovations is often in the yellow zone.

We conclude that the use of the $t_5$ distribution for the innovations slightly improve the results, and the VaR quality is acceptable for both Gaussian and scaled student quantiles, while the application of the empirical distribution of the residuals leads to an almost perfect fit of the prescribed quantiles for all considered time horizons.

Weierstrass-Institute
and
Humboldt University Berlin
Mohrenstr. 39
10117 Berlin
Germany
E-mail: spokoiny@wias-berlin.de